\documentclass[11pt,a4paper,reqno]{amsart}
\usepackage[T1]{fontenc}
\usepackage{lmodern}
\usepackage[hmargin=24mm,vmargin=21mm,headheight=14pt]{geometry}
\usepackage{amsmath,amssymb,mathtools,microtype}
\usepackage[unicode,colorlinks=true,linkcolor=black,citecolor=black,urlcolor=black,pdftitle={Boundaries and deformations of polyhedral Alexandrov spaces},pdfauthor={Nina Lebedeva and Polina Vorobyeva}]{hyperref}
\allowdisplaybreaks[2]
\newtheorem{theorem}{Theorem}[section]
\newtheorem{lemma}[theorem]{Lemma}
\newtheorem{proposition}[theorem]{Proposition}
\newtheorem{corollary}[theorem]{Corollary}
\theoremstyle{remark}

\DeclareMathOperator{\diam}{diam}
\DeclareMathOperator{\length}{length}
\DeclareMathOperator{\curv}{curv}
\DeclareMathOperator{\Int}{Int}

\DeclareMathOperator{\relint}{relint}
\newcommand{\CBB}{\mathrm{CBB}}
\newcommand{\intr}{\mathrm{intr}}
\newcommand{\R}{\mathbb R}
\newcommand{\Sph}{\mathbb S}

\numberwithin{equation}{section}
\title[Boundaries and deformations of polyhedral spaces]{Boundaries and deformations of\\polyhedral Alexandrov spaces}
\author{Nina Lebedeva}
\author{Polina Vorobyeva}
\date{November 30, 2026}
\begin{document}
\begin{abstract}
We prove that each boundary component of an $M_\kappa$-polyhedral Alexandrov space with curvature bounded below by $\kappa$ inherits the same lower curvature bound in its intrinsic length metric.
For polyhedral spaces of nonnegative curvature with nonempty boundary, we construct deformations that preserve nonnegative curvature and the intrinsic metric of the boundary. In the compact case, in dimension at least two, the resulting spaces are pairwise nonsimilar. The construction also applies to arbitrarily small strongly convex polyhedral neighborhoods of any point.
\end{abstract}
\maketitle

\section{Introduction}

The boundary conjecture is one of the major open problems in Alexandrov geometry. It asserts that the boundary of an Alexandrov space with curvature bounded below by $\kappa$, equipped with its intrinsic length metric, has curvature bounded below by $\kappa$; see~\cite[\S~13.4]{BGP} and~\cite[Section~9]{Petrunin}.
In general, the intrinsic metric of the boundary differs from the restriction of the ambient metric.

The conjecture has been proved in several special cases.
For convex hypersurfaces in Riemannian manifolds with sectional curvature bounded below by $\kappa$, the desired bound was obtained by Alexander, Kapovitch, and Petrunin~\cite{AKPconvex}, refining an earlier result of Buyalo.
Grove, Moreno, and Petersen~\cite{GMP} proved the conjecture for leaf spaces of singular Riemannian foliations with closed leaves on closed Riemannian manifolds. Their result includes orbit spaces of isometric actions of compact groups on such manifolds.
In~\cite{LebedevaSmoothable}, the corresponding bound was obtained for the boundaries of compact full-dimensional convex bodies in smoothable Alexandrov spaces (noncollapsing limits of Riemannian manifolds with the same lower sectional curvature bound).
In dimension two, the question reduces to the one-dimensional geometry of the boundary; for positive lower curvature bounds, an estimate on its length is also important.

In this paper we prove the conjecture for polyhedral spaces.
By an $M_\kappa$-polyhedral space we mean a complete length space admitting a locally finite triangulation in which each simplex, with the induced metric, is isometric to a simplex in the model space $M_\kappa$ of constant curvature $\kappa$; cf.~\cite[Definition~12.1]{AKP}.
We write $\partial P$ for the boundary of $P$ in the Alexandrov sense.

\begin{theorem}[Boundary theorem]\label{thm:main}
Let $P$ be an $n$-dimensional $M_\kappa$-polyhedral Alexandrov space, $n\ge2$, with $\curv P\ge\kappa$.
Then each component of $\partial P$, equipped with its intrinsic length metric, is complete and has curvature bounded below by $\kappa$.
\end{theorem}

For simplicity, we restrict the deformation results to the case $\kappa=0$.
For $\kappa\ne0$, one may consider analogous constructions using the spherical and hyperbolic models.

As an application of the boundary theorem, we construct deformations of polyhedral spaces with nonempty boundary that preserve nonnegative curvature and the intrinsic metric of the boundary.
The original space need not be compact.

\begin{theorem}[Deformation]\label{thm:graph}
Let $(P,d)$ be a Euclidean polyhedral Alexandrov space of nonnegative curvature with nonempty boundary.
There is a family of metrics $d_\varepsilon$ on $P$, $\varepsilon\ge0$, with $d_0=d$, such that each $(P,d_\varepsilon)$ is a Euclidean polyhedral space of nonnegative curvature.
The boundary remains the same subset of $P$, its intrinsic length metric is preserved, and for all $x,y\in P$,
\[
 d(x,y)\le d_\varepsilon(x,y)
       \le\sqrt{1+\varepsilon^2}\,d(x,y).
\]
If $P$ is compact and $\dim P\ge2$, the deformation gives a family of pairwise nonsimilar spaces.
\end{theorem}

Let us explain the difficulty in deforming polyhedral metrics while preserving a lower curvature bound.
Consider first a Euclidean polyhedral metric on the two-sphere and a triangulation in which the total angle at every vertex is strictly less than $2\pi$.
If the edge lengths are changed sufficiently little, the triangles remain nondegenerate and the angle sums at the vertices remain less than $2\pi$.
It is thus easy to construct a deformation preserving positive angle defects.
This argument relies on the strictness of the angle inequalities.
When the total angle is $2\pi$, as on a flat torus, even a small change of edge lengths may introduce negative curvature.
In dimension three, the same difficulty arises at edges around which the sum of dihedral angles is $2\pi$.
For spaces with boundary, Theorem~\ref{thm:graph} provides a deformation preserving nonnegative curvature without requiring strict angle inequalities.

The construction also applies locally. Every point of a Euclidean polyhedral $\CBB(0)$ space has arbitrarily small compact strongly convex polyhedral neighborhoods to which the deformation theorem can be applied.
The ambient space need not have boundary.
Here and throughout the paper, a subset is called strongly convex if it contains every minimizing geodesic between any two of its points.

The proof of the deformation theorem uses the boundary theorem (\S\,\ref{sec:boundary}) and the piecewise affine nature of the distance to the boundary (\S\,\ref{sec:PL}).
In \S\,\ref{sec:graph} we use Petrunin's construction involving graphs of concave functions~\cite[Section~9]{Petrunin}, and in \S\,\ref{sec:local} we construct small convex polyhedral neighborhoods using Voronoi cells.

\section{The boundary theorem}\label{sec:boundary}

The boundary of an Alexandrov space is a closed extremal subset.
By a theorem of Perelman and Petrunin~\cite[3.2(3)]{PerelmanPetrunin} (see also~\cite[\S~4.1]{Petrunin}), its intrinsic metric is locally bi-Lipschitz equivalent to the induced metric.
It follows that the boundary components are complete and locally compact in their intrinsic length metrics.
In the polyhedral case, the boundary consists of the free codimension-one faces and their subfaces, as follows from the definition in terms of spaces of directions.

\begin{lemma}\label{lem:connected}
The boundary of a compact spherical polyhedral space $X\in\CBB(1)$ of dimension at least two is connected if nonempty.
\end{lemma}
\begin{proof}
Suppose that $A$ and $B$ are distinct components of the boundary of $X$, and let $r=d(A,B)>0$.
Glue copies $X_i$, $i\in\mathbb Z$, along $A$ at the even-numbered seams and along $B$ at the odd-numbered seams, identifying corresponding points.
Each copy isometrically embeds into the resulting chain.
Indeed, transferring a finite chain of segments in different copies back to $X$ gives a path of the same length. By the definition of the gluing metric, leaving a copy cannot shorten the distance between two of its points.
A path crossing $m$ consecutive seams has length at least $(m-1)r$, since between successive seams it joins $A$ to $B$ in one of the copies.
Thus the chain is unbounded, and every bounded set lies in a finite union of compact copies. In particular, the resulting space is complete.

A sufficiently small neighborhood of each seam point is isometric to a neighborhood in the double of $X$, since the other boundary component is a positive distance away. Away from the seams the local metric agrees with that of $X$.
By Perelman's doubling theorem~\cite{BGP,Petrunin} and globalization, the resulting space is $\CBB(1)$.
Its unboundedness contradicts the diameter bound $\pi$ in dimension at least two~\cite[Theorems~8.31,8.44]{AKP}.
\end{proof}

\begin{lemma}\label{lem:links}
Let $P$ be an $M_\kappa$-polyhedral Alexandrov space, $\dim P\ge3$, and let $p\in\partial P$.
For sufficiently small $r>0$, the set $B_P(p,r)\cap\partial P$ is strongly convex in the intrinsic metric of the boundary and is isometric to the ball of radius $r$ centered at the tip of the $\kappa$-cone over $\partial\Sigma_pP$, equipped with its intrinsic length metric.
\end{lemma}
\begin{proof}
A sufficiently small ball $B_P(p,\rho)$ is isometric to a ball centered at the tip of the $\kappa$-cone over $\Sigma_pP$; see~\cite[Section~12A, p.~178]{AKP}.
In this chart, the boundary consists of cones over the faces of $\partial\Sigma_pP$.
The latter is connected by Lemma~\ref{lem:connected}.
Lengths of curves in these faces agree with their lengths in the $\kappa$-cone over $\partial\Sigma_pP$ with its intrinsic metric.

Let $r<\rho/4$.
Any two points of $B_P(p,r)\cap\partial P$ can be joined in the boundary by a path through $p$ of length less than $2r$; any path leaving $B_P(p,\rho)$ has length at least $2(\rho-r)$.
Hence it suffices to consider paths in the conical chart when computing intrinsic distances.
In a small $\kappa$-cone, radial projection onto a closed ball does not increase lengths and strictly shortens any rectifiable path with endpoints in the ball that leaves it: the radial contribution to length decreases and the angular contribution does not increase (for $\kappa>0$, we take $\rho<\pi/(2\sqrt\kappa)$).
Applying this to the ball whose radius is the larger of the distances from the two endpoints to $p$, we obtain both the desired isometry and strong convexity.
\end{proof}

\begin{proposition}[Perimeter bound]\label{prop:perimeter}
Let $Q$ be a compact connected two-dimensional spherical polyhedral Alexandrov space of curvature at least $1$ with nonempty boundary.
Then $Q$ is a disk and $\length(\partial Q)\le2\pi$.
\end{proposition}
\begin{proof}
Write $\alpha_i\le2\pi$ for the total angles at the interior vertices and $\beta_j\le\pi$ for the boundary angles.
By the Gauss--Bonnet formula,
\[
 2\pi\chi(Q)=\operatorname{area}(Q)
       +\sum_i(2\pi-\alpha_i)+\sum_j(\pi-\beta_j)>0.
\]
The classification of compact surfaces with boundary implies that $Q$ is a disk.
There are finitely many interior cone points with total angle strictly less than $2\pi$. We will eliminate them one by one without decreasing the perimeter.

First suppose there are no interior cone points.
On the simply connected open disk there is a spherical developing map $\operatorname{dev}:\Int Q\to\Sph^2$.
A minimizing geodesic between two interior points cannot have a first point of contact with the boundary: at such a point its opposite directions would have to be at angular distance $\pi$ in an interval of length $\beta\le\pi$.
This is possible only when $\beta=\pi$ and both directions are tangent to the boundary.
When $\beta=\pi$, the boundary is locally geodesic, so a minimizing geodesic tangent to it locally coincides with the boundary, contrary to the choice of the first point of contact.

By the diameter theorem, $\diam Q\le\pi$.
Under development, every minimizing geodesic between interior points becomes a minimizing great-circle arc of the same length.
The developing map is therefore distance preserving and, by completeness, extends to an isometric embedding of $Q$ as a convex spherical polygon.
By the spherical Crofton formula~\cite[Chapter~18]{Santalo}, its perimeter is at most $2\pi$, since almost every great circle meets the boundary of a convex polygon in at most two points.

In the general case, consider a shortest segment from an interior point $a$ to the boundary, and set $r=d_Q(a,\partial Q)>0$.
In the double, $d_{D(Q)}(a^+,a^-)=2r$: every path joining the two points meets the boundary, and the two copies of a shortest segment to the boundary together have length $2r$.
The doubling theorem and the diameter bound yield $r\le\pi/2$.

Now let $a$ be a cone point and let $[ax]$ be a shortest segment from $a$ to the boundary.
The open segment meets neither the boundary nor another cone point: at a point with total angle $\alpha'<2\pi$, the circle of directions has diameter $\alpha'/2<\pi$.
First variation in the two boundary directions at $x$ shows that both angles with $[ax]$ are at least $\pi/2$.
Their sum is the boundary angle, which is at most $\pi$.
Thus $x$ is a regular boundary point and $[ax]$ is perpendicular to the boundary. This remains true if $x$ is an artificial vertex of the triangulation.

Let $\alpha$ be the total angle at $a$.
Choose $0<\theta<\pi$ so that $\alpha+\theta\le2\pi$.
Cut the disk along $[ax]$ and insert an isosceles spherical triangle with side lengths $r$ adjacent to its apex angle $\theta$.
The triangle is nondegenerate. Since the spherical ball of radius $r\le\pi/2$ is convex, the base of the triangle lies in the ball.
At the endpoints of the base, the lateral sides point along the inward normals to the boundary of the ball; hence the base angles $\gamma$ are at most $\pi/2$.
Let $b$ denote the length of the base.

Along the glued banks of the cut the total angle is $2\pi$; at $a$ it is $\alpha+\theta$; at each new copy of $x$ the boundary angle is $\pi/2+\gamma\le\pi$.
All other angles are unchanged. The resulting spherical cone surface is again a compact disk. At every point the local spherical cone model has curvature at least $1$, by the angle bounds and the cone theorem~\cite[Theorem~11.6(a)]{AKP}. Globalization shows that the new disk is $\CBB(1)$.
Its perimeter has increased by exactly $b$.

To eliminate the cone point $a$, set $\theta=(2\pi-\alpha)/2\in(0,\pi)$.
After the first insertion the total angle at $a$ is $\pi+\alpha/2<2\pi$.
In the resulting $\CBB(1)$ disk, take a new shortest segment from $a$ to the boundary and repeat the insertion with the same $\theta$.
The total angle at $a$ becomes $2\pi$; the other interior cone points are unchanged and no new ones appear.
Eliminating all interior cone points in this way reduces the problem to the case already considered and proves the perimeter bound.
\end{proof}

The perimeter bound is also a special case of the general boundary volume estimate in~\cite[3.3.5]{Petrunin}.

\begin{proof}[Proof of Theorem~\ref{thm:main}]
There is nothing to prove if the boundary is empty. Normalize $\kappa$ to $-1,0,1$ and proceed by induction on the dimension, simultaneously for all three signs.
In dimension two, each boundary component with its intrinsic metric is a complete one-dimensional manifold without boundary, hence a line or a circle.
For $\kappa\le0$, it has the desired lower curvature bound.
For $\kappa=1$, the space $P$ is bounded by the diameter theorem and compact by completeness and local compactness. Its locally finite triangulation is therefore finite.
By Proposition~\ref{prop:perimeter}, its boundary is a single circle of length at most $2\pi$.
Under the convention in~\cite[Section~8H]{AKP}, the one-dimensional $\CBB(1)$ condition also allows longer circles, so the perimeter bound is used separately.

Now suppose $n\ge3$ and let $p\in\partial P$.
The compact spherical polyhedral space $\Sigma_pP$ has dimension $n-1$, curvature at least $1$, and nonempty connected boundary by Lemma~\ref{lem:connected}.
By the induction hypothesis, $L=(\partial\Sigma_pP,d_{\intr})$ is a complete $\CBB(1)$ space.
Moreover, $\diam L\le\pi$: for $n=3$ this follows from the perimeter bound, and for $n\ge4$ from the diameter theorem in dimension at least two.
By the cone theorem~\cite[Theorem~11.6(a)]{AKP}, the $\kappa$-cone over $L$ has curvature at least $\kappa$ near its tip.
Lemma~\ref{lem:links} identifies this metric with the intrinsic metric of the boundary of $P$.
The boundary components are complete, and globalization~\cite[Theorem~8.31]{AKP} completes the proof.
\end{proof}

\section{Distance to the boundary}\label{sec:PL}

In this and the following sections we consider Euclidean polyhedral spaces of nonnegative curvature.
We will prove that the distance to the boundary is piecewise affine.
The idea is to unfold the simplices along a shortest segment to the boundary, expressing its length as the distance to a hyperplane.
For each simplex, there are only finitely many possible unfoldings. Thus, at every point the distance function agrees with one of finitely many affine functions.
The following lemma makes this argument precise.

\begin{lemma}\label{lem:PL}
Let $P$ be a Euclidean polyhedral Alexandrov space of nonnegative curvature with nonempty boundary.
Then the function $f=d_P(\cdot,\partial P)$ is affine on the cells of some locally finite polyhedral subdivision. If $P$ is compact, the subdivision can be chosen finite.
\end{lemma}

We need two observations about minimizing geodesics for its proof.
We call a point flat if it has a neighborhood isometric to a Euclidean ball, and a flat boundary point if it has a neighborhood isometric to a Euclidean half-ball.

\begin{lemma}\label{lem:flat}
If one endpoint of a minimizing geodesic is flat, then every point in its interior is flat.
The interior of any minimizing geodesic either lies entirely in the boundary or is disjoint from it.
If $x$ is flat and $y$ is a nearest boundary point to $x$, then $y$ is a flat boundary point and $[xy]$ is perpendicular to the boundary.
\end{lemma}
\begin{proof}
At each interior point of a minimizing geodesic, its opposite directions determine a line in the tangent cone.
By the splitting theorem~\cite[Theorem~16.22]{AKP}, this cone has the form $\R\times C'$, with the line as its $\R$-axis.
In a conical neighborhood, translations along the axis identify small neighborhoods of nearby points of the geodesic.
Consequently, both flatness and membership in the boundary are locally constant along the open segment.
If one endpoint is flat, points of the segment sufficiently close to it are flat as well.

Consider the double of $P$. It is again Euclidean polyhedral: the copies are isometrically embedded and their triangulations agree along the boundary.
By the doubling theorem it has nonnegative curvature.
The union of $[xy]$ and its reflection is a minimizing geodesic joining the two copies of $x$, since every path joining them meets the boundary.
By what we have just proved, $y$ has a Euclidean neighborhood in the double.
The involution interchanging the copies acts in this neighborhood as an orthogonal involution whose fixed-point set separates the ball into two parts.
It is therefore a reflection in a hyperplane; $y$ is a flat boundary point, and the doubled segment is perpendicular to this hyperplane.
\end{proof}

\begin{lemma}[Unfoldings]\label{lem:galleries}
For minimizing geodesics with flat interiors contained in a fixed compact set, after barycentric subdivision a finite list of sequences of adjacent full-dimensional simplices suffices for unfolding along these geodesics.
At a flat endpoint or a flat boundary endpoint, the unfolding can be extended to any full-dimensional cell containing that endpoint.
\end{lemma}
\begin{proof}
Every simplex of the barycentric subdivision is strongly convex.
Indeed, it corresponds to a chain of nested faces of the original triangulation.
For two points $u,v$ in such a simplex, the smallest original faces containing them are nested. We may assume that $u\in F$ and $v\in\relint F$.
Extend the chord $uv$ in $F$ slightly past $v$ to a point $w$.
Since $F$ is isometrically embedded, $d(u,w)=|uv|+|vw|$.
Any other minimizing geodesic from $u$ to $v$, extended by $vw$, would give two minimizing geodesics starting at $w$ that initially coincide and then diverge, contrary to nonbranching~\cite[Theorem~8.37]{AKP}.

Thus a minimizing geodesic meets each simplex in a single (possibly degenerate) segment.
Only finitely many simplices meet a compact set.
When a geodesic passes through a face of codimension greater than one, we complete the sequence of simplices by a chain of full-dimensional simplices in the star of the point of passage, consecutive simplices sharing a facet.
Such a chain exists because the star gives a subdivision of a Euclidean ball near a flat point: a polygonal path in general position crosses only facets.
The same argument applies in a half-ball at a flat boundary endpoint.
Choosing chains without repetitions gives a uniform bound on the number of simplices in an unfolding.
All of them meet the original compact set. Hence there are only finitely many possible sequences.
The unfoldings agree with the flat local charts, so the geodesic unfolds to a straight segment.
\end{proof}

\begin{proof}[Proof of Lemma~\ref{lem:PL}]
Fix a full-dimensional simplex $S$ of the barycentric subdivision.
Minimizing geodesics from points of $S$ to the boundary have length at most $\max_S f$ and therefore lie in a common compact set (a complete locally compact length space is proper; see~\cite{BBI}).
If $x\in\Int S$, then by Lemma~\ref{lem:flat} a shortest segment to the boundary has flat interior and meets the boundary perpendicularly at a flat boundary point.
By Lemma~\ref{lem:galleries}, only finitely many boundary hyperplanes arise when unfolding into the chart of $S$.
Signed distances to these hyperplanes, with both choices of normal, give affine functions $\ell_1,\ldots,\ell_N$ such that
\[
 f(x)\in\{\ell_1(x),\ldots,\ell_N(x)\}.
\]
Cut $S$ along the hyperplanes where pairs of these functions agree.
By continuity, in the interior of each resulting full-dimensional cell, $f$ agrees with a single function $\ell_i$; the equality extends to the closure of the cell.
On common faces of simplices, take common refinements of the subdivisions, starting with the lower-dimensional faces, and extend them to the cells.
Local finiteness is preserved.
In particular, $|\nabla f|=1$ on full-dimensional affine cells.
\end{proof}

\section{Deformation}\label{sec:graph}

Let $P$ satisfy the assumptions of Theorem~\ref{thm:graph}, and set $f=d_P(\cdot,\partial P)$.
By the theorem on distance to the boundary~\cite[Theorem~3.3.1]{Petrunin}, $f$ is concave in the interior of $P$.
By Lemma~\ref{lem:flat}, the interior of a minimizing geodesic either lies in the boundary, where $f=0$, or lies entirely in the interior of $P$.
By continuity, $f$ is concave along all minimizing geodesics in $P$.

\begin{proof}[Proof of Theorem~\ref{thm:graph}]
For $\varepsilon>0$, consider the set
\[
 W_\varepsilon=\{(x,t)\in P\times\R:|t|\le\varepsilon f(x)\}.
\]
By concavity of $f$, it is closed and strongly convex in the complete $\CBB(0)$ space $P\times\R$~\cite[Theorem~11.6(b)]{AKP}; hence, with its intrinsic metric, it is complete and has nonnegative curvature.
By Lemma~\ref{lem:PL}, the part of $W_\varepsilon$ over each affine cell is a Euclidean polyhedron.
A compatible triangulation of these polyhedra gives a locally finite metric triangulation of $W_\varepsilon$: distances in the cells are Euclidean because the base cells are isometrically embedded and $W_\varepsilon$ is convex.
Its free facets form two graphs
$\Gamma_\varepsilon^{\pm}=\{(x,\pm\varepsilon f(x)):x\in P\}$, glued along $\partial P\times\{0\}$.

By Theorem~\ref{thm:main}, the boundary of $W_\varepsilon$ has nonnegative curvature in its intrinsic metric.
The upper graph $\Gamma_\varepsilon=\Gamma_\varepsilon^+$, with its own intrinsic metric, embeds isometrically into this boundary: reflecting portions of a path from the lower graph to the upper graph gives a path of no greater length.
Thus $\Gamma_\varepsilon$ also has nonnegative curvature, by four-point comparison.

Transfer the metric of the graph to $P$ via $i_\varepsilon(x)=(x,\varepsilon f(x))$.
Projection from the graph to $P$ does not increase lengths, and lifting any rectifiable path in $P$ increases its length by a factor of at most $\sqrt{1+\varepsilon^2}$, since $f$ is $1$-Lipschitz.
It follows that
\begin{equation}\label{eq:distances}
 d(x,y)\le d_\varepsilon(x,y)\le\sqrt{1+\varepsilon^2}\,d(x,y).
\end{equation}
In particular, the graph is complete.
The graphs over the affine cells form a polyhedral complex with the same combinatorial structure as the subdivision of $P$.
These cells are isometrically embedded: the distance in the product between two points of such a cell equals the length of their Euclidean chord, and the chord itself lies in the graph.
By the description of the boundary in terms of free facets, $i_\varepsilon$ identifies the boundaries of $P$ and the graph; their intrinsic metrics agree because $f=0$ on the boundary.
For $\varepsilon=0$, we recover the original metric.

On a full-dimensional affine cell, lengths in the gradient direction of $f$ are multiplied by $\sqrt{1+\varepsilon^2}$, while lengths in orthogonal directions are preserved.
Thus, for $n=\dim P$ and any Borel set $E\subset P$,
\begin{equation}\label{eq:volume}
 \mathcal H^n(i_\varepsilon(E))
   =\sqrt{1+\varepsilon^2}\,\mathcal H^n(E).
\end{equation}
Here we used $|\nabla f|=1$ and the fact that the lower-dimensional skeleton has zero $n$-dimensional measure.
If $P$ is compact and $n\ge2$, then its boundary has finite positive $(n-1)$-dimensional measure.
Since this measure is preserved, any similarity between two of the resulting spaces must have scale factor one.
But by~\eqref{eq:volume}, their volumes are different for distinct parameters.
\end{proof}

In the noncompact case, the spaces may also be pairwise nonsimilar.
For example, if $P=\R_+^2$, then $f(x,y)=\min(x,y)$.
The graph consists of two flat sectors glued along the lifted diagonal.
Unfolding gives a sector of angle $2\arccos(1/\sqrt{2+\varepsilon^2})$, which is strictly increasing in $\varepsilon$.
For a half-space, on the other hand, the graph is a flat half-space isometric to the original one.

\section{Local deformation}\label{sec:local}

To apply the deformation theorem near an arbitrary point, we construct small strongly convex polyhedral neighborhoods.
Instead of balls, which are generally not polyhedral, we use Voronoi cells in the tangent cone.
This approach is related to the method of convex cells in~\cite[Sections~4--5]{LebedevaPetrunin}.

Let $C=T_pP$ be the tangent cone with tip $o$, and let $\Sigma=\Sigma_pP$.
Subdivide the triangulation to make $p$ a vertex. Extending the simplices in its finite star radially from $p$ gives a finite decomposition of $C$ into Euclidean simplicial cones.
These cones are isometrically embedded: any two points in one cone can be rescaled into a small metric simplex in the conical neighborhood of $p$.
To apply Lemma~\ref{lem:galleries}, subdivide these cones into bounded polyhedra.
On each cone, take the linear function $\tau$ equal to $1$ on its unit generators and cut along the hyperplanes $\tau=j$, $j=1,2,\ldots$.
The functions $\tau$ agree on common faces and are comparable to the distance from the tip.
A compatible triangulation gives a locally finite metric triangulation of $C$.

\begin{lemma}\label{lem:radial}
For a regular direction $v\in\Sigma$, the function $h_v(ru)=r\cos d_\Sigma(u,v)$ is linear on the cells of some finite conical subdivision of $C$.
\end{lemma}
\begin{proof}
Set $q=(1,v)$; the point $q$ is flat.
Fix a full-dimensional cone $S$ in the decomposition.
For a flat unit point $x\in S$, a minimizing geodesic to $q$ lies in the ball of radius $3$ and has flat interior by Lemma~\ref{lem:flat}.
By Lemma~\ref{lem:galleries}, there are only finitely many possible unfoldings into the chart of $S$.
The transitions between the charts of the original cones fix $o$, so each unfolding gives a unit vector $q_G$.
By the cone distance formula and the equality $d_C(x,q)=|x-q_G|$ for the unfolding in question,
\[
 h_v(x)=\tfrac12(|x|^2+1-d_C(x,q)^2)=\langle x,q_G\rangle.
\]
By continuity and homogeneity, at every point of $S$, the function $h_v$ agrees with one of finitely many linear functions.
As in the proof of Lemma~\ref{lem:PL}, cutting the cones along the hyperplanes where pairs of these functions agree and taking compatible subdivisions gives the required finite conical decomposition.
\end{proof}

\begin{proposition}\label{prop:neighborhood}
Every point of a Euclidean polyhedral Alexandrov space of nonnegative curvature has arbitrarily small compact strongly convex polyhedral neighborhoods.
\end{proposition}
\begin{proof}
The assertion is clear in dimension less than two. Suppose $\dim P\ge2$.
Choose a finite $\delta$-net of regular directions $v_i$ in $\Sigma$, where $0<\delta<\pi/2$.
For $R>0$, set $q_i=(2R,v_i)$ and consider the Voronoi cell of the tip $o$ with respect to the points $o,q_1,\ldots,q_m$:
\[
 D_R=\{x\in C:d_C(x,o)\le d_C(x,q_i)\text{ for all }i\}.
\]
By the cone distance formula,
\[
 D_R=\bigcap_i\{h_{v_i}\le R\}.
\]
The function $-h_v$ is the Busemann function of the ray $(t,v)$:
\[
 -h_v(x)=\lim_{t\to\infty}(d_C(x,(t,v))-t).
\]
In a space of nonnegative curvature it is concave~\cite[Exercise~8.25]{AKP}.
Hence $D_R$ is strongly convex.
By the net property, $r\cos\delta\le\max_i h_{v_i}(ru)\le r$. Thus
\begin{equation}\label{eq:balls}
 \overline B_R(o)\subset D_R\subset\overline B_{R/\cos\delta}(o),
 \qquad D_R=R D_1.
\end{equation}
By Lemma~\ref{lem:radial}, we can subdivide $C$ into finitely many metric cones on each of which all the $h_{v_i}$ are linear.
The inequalities $h_{v_i}\le R$ cut out bounded Euclidean polyhedra in these cones.
Their compatible triangulation is metric, since by strong convexity the intrinsic metric of $D_R$ agrees with the restriction of the cone metric.
The neighborhood $D_R$ has nonempty boundary, since in the regular part of the cone there are exterior free facets of these polyhedra.

Transfer $D_R$ into the conical neighborhood $B_P(p,\rho)$.
If $3R/\cos\delta<\rho$, every minimizing geodesic between points of $D_R$ stays in this neighborhood: its length is at most $2R/\cos\delta$, by taking the path through $p$.
In the conical chart it lies in $D_R$ by strong convexity.
The inclusions~\eqref{eq:balls} show that as $R\to0$, these sets form a neighborhood basis at $p$.
\end{proof}

\begin{corollary}[Local deformation]\label{cor:local}
Every point of a Euclidean polyhedral Alexandrov space of nonnegative curvature and dimension at least two has arbitrarily small compact strongly convex polyhedral neighborhoods admitting a family of pairwise nonsimilar polyhedral metrics of nonnegative curvature with the same intrinsic boundary metric. These metrics converge to the original one in the bi-Lipschitz sense.
\end{corollary}
\begin{proof}
The neighborhood $D_R$ from Proposition~\ref{prop:neighborhood} has nonnegative curvature by strong convexity.
Theorem~\ref{thm:graph} applies to it.
Moreover, since $D_R=RD_1$, the distances to their respective boundaries satisfy $f_{D_R}(Rx)=R f_{D_1}(x)$.
Therefore the resulting metrics satisfy $d_{\varepsilon,R}(Rx,Ry)=R\,d_{\varepsilon,1}(x,y)$.
Thus for sufficiently small $R$ we obtain the same deformation up to rescaling.
\end{proof}

\medskip
The proof of the boundary theorem, including the geometric argument for the perimeter bound, builds on P.~Vorobyeva's bachelor's thesis~\cite{Vorobyeva}.

\medskip
\noindent\textbf{Use of AI.} 
The authors used ChatGPT (OpenAI) to assist with translating, drafting, and editing the manuscript, 
as well as checking the proofs for clarity and accuracy. All mathematical results and arguments were independently verified by the authors, who take full responsibility for the final text and content.

\end{document}